\documentclass[11pt,oneside,reqno]{amsart}
\usepackage[all]{xy}
\usepackage{amsfonts,amsmath,oldgerm,amssymb,amscd, comment,multirow}
\UseComputerModernTips
\numberwithin{equation}{section}
\usepackage[breaklinks]{hyperref}
\allowdisplaybreaks

\usepackage{enumerate}

\usepackage{caption} 
\newcommand{\mex}{\operatorname{mex}}

\newtheorem{example}{Example}
\newtheorem{theorem}{Theorem}[section]
\newtheorem{prop}{Proposition}

\newtheorem{lemma}[theorem]{{\bf Lemma}}

\begin{document}
	\vspace{-4cm}	
	\title[ The Minimal excludant of $\ell$-regular partition ]{The Minimal excludant of $\ell$-regular partition} 
		\author[Nabin Kumar Meher*]{Nabin Kumar Meher}
	\address{Nabin Kumar Meher, Department of Mathematics, National Institute of Technology Raipur, G. E. Road, Raipur, Chhatishgarh, India 492010.}
	\email{mehernabin@gmail.com, nkmeher.maths@nitrr.ac.in}
	\maketitle
	\pagenumbering{arabic}
	\pagestyle{headings}
\begin{abstract}
	Andrews and Newman \cite{andrewsnewmanI} introduced and studied ``mex,'' function for an integer partition $\lambda$ of a positive integer $n$. The related function $\sigma^{\mathrm{mex}}(n)$ was subsequently defined, and its generating function and other properties have been studied by several authors. In this article, for a positive integer $\ell \geq 2$, we define the minimal excludant of an $\ell$-regular partition $\lambda$ of a positive integer $n$. We then define $\sigma_\ell^{\mathrm{mex}}(n)$, the sum of the minimal excludants over all $\ell$-regular partitions of $n$, and obtain explicit generating functions for $\sigma_\ell^{\mathrm{mex}}(n)$ in the cases $\ell = 2$ and $\ell = 3$, expressed as eta-quotients. We prove that the sequence $\sigma_\ell^{\mathrm{mex}}(n)$ is weakly increasing for every $\ell$, and we obtain Hardy--Ramanujan type asymptotic formulas for $\sigma_2^{\mathrm{mex}}(n)$ and $\sigma_3^{\mathrm{mex}}(n)$ as $n \to \infty$. Additionally, we completely determine the parity of $\sigma_2^{\mathrm{mex}}(n)$, showing that it is odd precisely at generalized pentagonal numbers, and we use this characterization to compute the density of $n$ for which $\sigma_2^{\mathrm{mex}}(n)$ is odd versus even. Finally, for $\sigma_3^{\mathrm{mex}}(n)$, we prove that it is odd infinitely often, and we establish a quantitative lower bound, showing that the number of such $n$ up to $X$ grows at least like $\log \log X$.
\end{abstract}

	\maketitle
	\section{Introduction and Statement of Results}
	
	A partition $\pi = (\pi_1, \pi_2, \dots, \pi_k)$ of a non-negative integer $n$ is a non-increasing sequence of positive integers summing to $n$. Let $p(n)$ denote the number of partitions of $n$, and let $\mathcal{P}(n)$ denote the set of all partitions of $n$. Grabner and Knopfmacher \cite{grabner} introduced and studied an interesting partition statistic which they called the \emph{smallest gap} of a partition, namely the least positive integer missing from its parts. 
	
	In 2019, Andrews and Newman \cite{andrewsnewmanI, andrewsnewmanII} revisited this idea, thereby rediscovering the Grabner--Knopfmacher statistic, and went on to generalize it very naturally to arithmetic progressions. They defined
	\[
	\sigma\mathrm{mex}(n) := \sum_{\pi \in \mathcal{P}(n)} \mathrm{mex}(\pi),
	\]
	and proved the striking identity 
	\begin{align}\label{eq1.1}
	\sigma\mathrm{mex}(n) = D_2(n).
\end{align}
	 where $D_2(n)$ denotes the number of two-colored partitions of $n$ into distinct parts; the corresponding generating function had, in fact, already appeared implicitly in the work of Grabner and Knopfmacher. They also obtained a Hardy--Ramanujan--Rademacher type exact formula for $\sigma\mathrm{mex}(n)$ in terms of Bessel functions and Dedekind sums, together with the asymptotic estimate $\sigma\mathrm{mex}(n) \sim \tfrac{1}{4\sqrt{6n^3}}\exp\!\big(\pi\sqrt{2n/3}\big)$ as $n \to \infty$. In the same paper, they introduced a second arithmetic function
	\[
	a(n) := \sum_{\substack{\pi \in \mathcal{P}(n) \\ \mathrm{mex}(\pi) \text{ odd}}} 1,
	\]
	and showed that $\sigma\mathrm{mex}(n) \equiv a(n) \pmod 2$, with $a(n)$ is almost always even and is odd exactly when $n$ is of the form $r(3r\pm 1),$ for some positive integer $r$.
	In 2021, Ballantine and Merca \cite{ballantine-merca} gave a combinatorial proof of the Andrews--Newman identity \eqref{eq1.1}. Barman and Singh \cite{barman-singh1, barman-singh2}, and later da Silva and Sellers \cite{sellers-dasilva}, established further congruences and parity results for mex-related partition functions. In 2021, Chern \cite{chernmaxexcludnt} introduced the notion of the \emph{maximal excludant} $\mathrm{maex}(\pi)$, the largest integer below the biggest part of $\pi$ that fails to occur in $\pi$, and obtained both a generating function identity for $\sigma\mathrm{maex}(n) := \sum_{\pi \in \mathcal{P}(n)} \mathrm{maex}(\pi)$ and the asymptotic relation $\sigma\mathrm{maex}(n) \sim \sigma L(n)$, where $\sigma L(n)$ is the sum of the largest parts over all partitions of $n$, with $\sigma L(n)$ itself satisfying a Hardy--Ramanujan type formula due to Kessler and Livingston \cite{Kessler-Livingston}.
	In 2022, Kaur, Bhoria, Eyyunni and Maji \cite{KaurBhoriaEyyunniMaji} have extended these investigations to mex statistics on partitions restricted in various ways, such as partitions into distinct parts and they studied the sum of minimal excludants over partitions into distinct parts $\sigma_d \mathrm{mex}(n)$ and found a connection of Ramanujan's sum $\sigma(q)$ with the sum function $\sigma_d \mathrm{mex}(n).$ In 2026, Baruah and Talukdar \cite{baruah-Talukdar} studied the partition functions involving least $r$ gaps and obtained the generating function of the related sum function in terms of eta quotients, studied their asymptotic formula.
	This growing literature makes clear that the minimal excludant, though a simple combinatorial rule, interacts richly with $q$-series, modular forms, and analytic number theory.
	 Motivated by the above results, we study the minimal excludant of $\ell$-regular partitions. 

 Let $\ell \geq 2$ be a positive integer. Let $P_{\ell}(n)$ denote the set of all $\ell$-regular partitions of $n$ and $ p_{\ell}(n) := \left|\mathcal{P}_{\ell}(n)\right|. $ For a partition $\lambda\in P_{\ell}(n)$, the
\emph{$\ell$-regular minimal excludant} is defined by
$\operatorname{mex}_{\ell}(\lambda)
:= \min\{\,r\ge1:\ \ell\nmid r,\ r\notin\lambda\,\}.$
In other words, $\operatorname{mex}_{\ell}(\lambda)$ is the smallest positive integer that is not divisible by $\ell$ and does not occur as a part of $\lambda$. Let $m$ be positive integer  such that $ \ell\nmid m.$  We define
\begin{align*}
	p_{\ell}^{\operatorname{mex}}(m,n)
	:=
	\#\left\{
	\lambda\in P_{\ell}(n):
	\operatorname{mex}_{\ell}(\lambda)=m
	\right\}.
\end{align*}  If $\ell\mid m$, then $ p_{\ell}^{\operatorname{mex}}(m,n)=0,$ by convention. We consider some examples to understand the above definitions.
\begin{example}
	Let $\ell=3$ and $n=12$. Then, the admissible parts are $ 1,2,4,5,7,8,10,11,\ldots.$ 
	\begin{enumerate}
		\item Consider $\lambda=(8,2,2).$ 
		The missing admissible parts are $ 1,\;4,\;5,\;7,\ldots.$
		Hence $\operatorname{mex}_{3}(\lambda)=1.$
		\item Let $ \lambda=(7,2,2,1).$ The missing admissible parts are $4,\; 5, \;8, \dots. $ Therefore $\operatorname{mex}_{3}(\lambda)=4.$
		\item Let$ \lambda=(5,4,2,1).$ The missing admissible integers are $7,\; 8, \; 10, \ldots.$	Hence $\operatorname{mex}_{3}(\lambda)=7.$ Notice that $\operatorname{mex}_{\ell}$ is unbounded.
	\end{enumerate}	
\end{example}

%
	
	Writing $\mathcal{P}_\ell(n)$ for the set of $\ell$-regular partitions of $n$, we let $p_\ell^{\mathrm{mex}}(m,n)$ denote the number of $\lambda \in \mathcal{P}_\ell(n)$ with $\mathrm{mex}_\ell(\lambda) = m, \quad \ell \nmid m$. In this article, we derive a general two-variable generating function for this statistic in  Proposition \ref{bivarategeneratingfunction}. Our main object of study is the aggregate function
	\[
	\sigma_\ell^{\mathrm{mex}}(n) := \sum_{\lambda \in \mathcal{P}_\ell(n)} \mathrm{mex}_\ell(\lambda),
	\]
	the natural analogue of $\sigma\mathrm{mex}(n)$ for $\ell$-regular partitions.
	 We obtain explicit generating functions for $\sigma_\ell^{\mathrm{mex}}(n)$ expressed as eta-quotients for $\ell = 2$ and $\ell = 3$. 
	 
 Let $\lambda\in \mathcal{P}_{2}(n).$ Suppose that $\operatorname{mex}_{2}(\lambda)=2r-1.$ Then $1,3,5,\ldots,2r-3$ all occur as parts of $\lambda$, whereas $2r-1$ does not occur.
	 We define
	 $$
	 p_{2}^{\widetilde{\operatorname{mex}}}(r,n)
	 :=
	 \#\left\{
	 \lambda\in P_{2}(n):
	 \operatorname{mex}_{2}(\lambda)=2r-1
	 \right\}.
	 $$ 
	 
	 \begin{theorem}\label{sigmamex2}
	 	We have
	 	\begin{align}
	 		\sum_{n\geq0} \sigma_{2}^{\mex}(n)q^n = \frac{(q^2;q^2)^6_{\infty}}{(q;q)^3_{\infty} (q^4;q^4)^2_{\infty}} = \frac{f_2^6}{f_1^3 f_4^2},
	 	\end{align}
	 	where, for brevity, we set $f_r= \left(q^r,q^r\right)_{\infty},$ $r \geq 1$ is a positive integer.
	 \end{theorem}
 	\begin{theorem}\label{sigmamex3}
 	We have
 	\begin{align}
 		\sum_{n\geq0} \sigma_{3}^{\mex}(n)q^n 
 		= \frac{f_6^5}{f_1 f_3 f_{12}^2} + q \frac{f_{12}^2f_3}{f_1 f_6}.
 	\end{align}
 \end{theorem}
	Further, we show that $\{\sigma_\ell^{\mathrm{mex}}(n)\}_{n \geq 0}$ is weakly increasing for every $\ell \geq 2.$ In particular, we find the following result.
	 	\begin{theorem}\label{monotonicitytheorem}
	 	The sequence $\{ \sigma_{\ell}^{\mex}(n) \}_{n\geq 0}$ is weakly increasing. Moreover, it satisfies the following inequality
	 	\begin{align}
	 		\sigma_{\ell}^{\mex}(n+1) - \sigma_{\ell}^{\mex}(n) \geq p_{\ell}(n)-p_{\ell}(n-1).
	 	\end{align}
	 \end{theorem}
	 Using Tauberian theorem of Ingham, we establish Hardy--Ramanujan type asymptotic formulas for $\sigma_2^{\mathrm{mex}}(n)$ and $\sigma_3^{\mathrm{mex}}(n).$ In particular, we have the following results.
	 \begin{theorem}\label{asymptotictheoremofsigma2}
	 	We have
	 	\begin{align*}
	 		\sigma_{2}^{\mex}(n)
	 		\sim
	 		\frac{1}{2\sqrt{2n}}\,
	 		e^{\pi\sqrt{n/3}},
	 		\qquad \hbox{as} \qquad n\to\infty.
	 	\end{align*}
	 \end{theorem}
 \begin{theorem}\label{asymptotictheoremofsigma3}
 	We have
 	\begin{align*}
 		\sigma_{3}^{\mex}(n)
 		\sim 
 		\frac{1}{4\sqrt{n}}\,
 		e^{\frac{2\pi}{3}\sqrt{n}},
 		\qquad \hbox{as} \qquad  n\to\infty.
 	\end{align*}
 \end{theorem}
We further determine the parity of $\sigma_2^{\mathrm{mex}}(n)$ completely, showing that it is odd precisely at generalized pentagonal numbers, exactly as in the Andrews--Newman setting, and deduce the corresponding density of $n$ for which $\sigma_2^{\mathrm{mex}}(n)$ is odd versus even. In particular, we get the following results.

\begin{theorem}\label{distributionofsigma2}
	For every real number \(X\geq 1\), we have
	\[
	\#\left\{1\leq n\leq X:\sigma_{2}^{\mex}(n) \ \text{is an odd integer}\right\}
	=
	\frac{2\sqrt{6X}}{3} +\mathcal{O}(1),
	\]
	while
	\[
	\#\left\{1\leq n\leq X:\sigma_{2}^{\mex}(n) \ \text{is an even integer}\right\}
	=
	X-\frac{2\sqrt{6X}}{3} +\mathcal{O}(1).
	\]	
\end{theorem}
Finally, we prove that $\sigma_3^{\mathrm{mex}}(n)$ is odd infinitely often, with a quantitative lower bound of order $\log \log X$ for the number of such $n$ up to $X$. In particular, we obtain the following result. 
		\begin{theorem}\label{distributionofsigma3} 
		For sufficiently large \(X\), we have
		\[
		\#\left\{
		0\leq n\leq X:
		\sigma_3^{\mathrm{mex}}(n)\equiv1\pmod{2}
		\right\}
		\gg\log\log X.
		\]
	\end{theorem} 
	
%

	\section{Preliminaries}
	\begin{lemma}\label{lem1}
		We have 
		\begin{align*}
			\sum_{t\ \in \mathbb{Z}} q^{t^2} = \frac{(q^2;q^2)^5_\infty }{ (q;q)^2_\infty (q^4;q^4)^2_\infty} = \frac{f_2^5}{f_1^2 f_4^2}.
		\end{align*}
	\end{lemma}
\begin{proof}
	See (\cite{Hirschhorn}Eq. ~(1.5.1)).
	\end{proof}
	\begin{lemma}\label{lemmaJacobitripleproduct}
		For $z \neq 0, |Q|<1,$ we have
		\begin{align}\label{tripleproduct}
			\sum_{r=-\infty}^{\infty} Q^{r^2} z^r = \left(-Qz; Q^2\right)_{\infty} \left(-\frac{Q}{z}; Q^2\right)_{\infty} \left(Q^2; Q^2\right)_{\infty}.
		\end{align}
	\end{lemma}
\begin{proof}
		See  (\cite{andrewspartitionbook} Theorem $2.8$).
\end{proof}
	\begin{lemma}[\cite{Ingham}]\label{LemmaIngham}
	Let
	$B(q):=\sum_{n=0}^{\infty}b(n)q^n$
	be a power series with radius of convergence $1$. Assume that $\{b(n)\}$ is a weakly increasing sequence of non-negative real numbers. If there are constants $\beta,\gamma\in\mathbb{R}$ and $\mu>0$ such that
	\[
	B(e^{-y})\sim \beta y^\gamma e^{\frac{\mu}{y}},
	\qquad \text{as }y\to0^+,
	\]
	then we have
	\[
	b(n)\sim
	\frac{\beta}{2\sqrt{\pi}}
	\frac{\mu^{\frac{2\gamma+1}{4}}}
	{n^{\frac{2\gamma+3}{4}}}
	e^{2\sqrt{\mu n}}
	\qquad \text{as }n\to\infty.
	\]
\end{lemma}
\begin{lemma}\label{lemmabinomial}
	For any positive integer $u$ and $j,$ and a prime number $p,$ we obtain
	\begin{align}
		\left(q^u;q^u\right)_{\infty}^{p^j} \equiv \left(q^{pu};q^{pu} \right)_{\infty}^{p^{j-1}} \pmod{p^j}.
	\end{align}
\end{lemma}
\begin{lemma}\label{pentagonaltheorems}
Euler's pentagonal theorem gives the expansion
	\begin{align}\label{Eulerpentagonaltheorem}
f_1 = \left(q;q\right)_{\infty} = \sum_{n=-\infty}^{\infty} (-1)^n q^{\frac{n(3n-1)}{2}} = 	1+\sum_{n=1}^{\infty} (-1)^n \left( q^{\frac{n(3n-1)}{2}} + q^{\frac{n(3n+1)}{2}} \right).
	\end{align}
	\end{lemma}
\begin{proof}
	See (\cite{andrewspartitionbook} Corollary $1.7$).
	\end{proof}
\begin{lemma}
	We have
	\begin{align}
f_1^3 = (q;q)_{\infty}^3=\sum_{n=0}^{\infty} (-1)^n(2n+1) q^{\frac{n(n+1)}{2}}.
\end{align}
\begin{proof}
	see (\cite{Hirschhorn}, ~Eq. (1.7.1) )
	\end{proof}
\end{lemma}
\section{Proofs}
	\begin{prop}\label{bivarategeneratingfunction}
		For a positive integer $\ell \geq 2,$
		define
		\[
		A_{\ell}(m):=\sum_{\substack{1\le j<m\\ \ell \nmid j}} j.
		\] Then we have
			\begin{align}\label{eq99a} 
			\sum_{\substack{m\ge1\\ \ell\nmid m}}
			\sum_{n=0}^{\infty}
			p_{\ell}^{\operatorname{mex}}(m,n)
			z^m q^n
			=
			\frac{(q^{\ell};q^{\ell})_{\infty}}
			{(q;q)_{\infty}}
			\sum_{\substack{m\ge1\\ \ell\nmid m}}
			z^m
			q^{A_{\ell}(m)}(1-q^m).
		\end{align}
		\end{prop}
	\begin{proof}
		 Suppose $ \operatorname{mex}_{\ell}(\lambda)=m, \quad \ell\nmid  m.$
		Then \begin{enumerate}
			\item every admissible integer $	1\le r<m,\quad \ell\nmid r,$ must occur at least once;
			\item the part $m$ is absent;
			\item every admissible integer greater than $m$ may occur with arbitrary
			multiplicity.
		\end{enumerate}
		The compulsory contribution of the parts below $m$ is
	\begin{align*}
		q^{A_{\ell}(m)} \prod_{\substack{r<m\\ \ell\nmid r}} \frac{1}{1-q^r}, \quad \hbox{and the unrestricted larger admissible parts contribute} \quad \prod_{\substack{r >m \\ \ell\nmid r}} \frac{1}{1-q^r}.
	\end{align*}
Thus, we obtain	
\begin{align}\label{eq98}
		\sum_{n=0}^{\infty} p_{\ell}^{\operatorname{mex}}(m,n) q^n
		= q^{A_{\ell}(m)} \prod_{\substack{r\geq 1\\ \ell\nmid r}}
		\frac{1-q^m}{1-q^r}, \quad
		\ell\nmid m.
	\end{align}
Since 
\begin{align}\label{eq98a}
\prod_{ \substack{r \geq 1 \\ \ell\nmid r}}\frac{1}{1-q^r}
=\frac{(q^{\ell};q^{\ell})_{\infty}}{(q;q)_{\infty}},	
\quad \hbox{we get} \quad 	\sum_{n=0}^{\infty} p_{\ell}^{\operatorname{mex}}(m,n) q^n
	= q^{A_{\ell}(m)} \frac{(q^{\ell};q^{\ell})_{\infty}}
	{(q;q)_{\infty}}(1-q^m), \quad
	\ell\nmid m.
	\end{align}
Consequently, we have
	\begin{align}\label{eq99} 
		\sum_{\substack{m\ge1\\ \ell\nmid m}}
		\sum_{n=0}^{\infty}
		p_{\ell}^{\operatorname{mex}}(m,n)
		z^m q^n
		=
		\frac{(q^{\ell};q^{\ell})_{\infty}}
		{(q;q)_{\infty}}
		\sum_{\substack{m\ge1\\ \ell\nmid m}}
		z^m
		q^{A_{\ell}(m)}(1-q^m).
	\end{align}
	\end{proof}
	
	\subsection{Proof of Theorem \ref{sigmamex2} and \ref{sigmamex3} }

\begin{proof}[Proof of Theorem \ref{sigmamex2}]
	Let $ p_{2}^{\widetilde{\operatorname{mex}}}(r,n) $ be the number of $2$-regular partitions $\lambda$ of $n$ such that $\operatorname{mex} (\lambda)=2r-1.$ We have
	\begin{align}\label{eq099}
		 G_{2}(z,q) := \sum_{n\geq0}\sum_{r\geq1} p_{2}^{\widetilde{\operatorname{mex}}}(r,n) z^{2r-1}q^n = \sum_{r=1}^{\infty} z^{2r-1} \sum_{n=0}^{\infty}
		 p_{2}^{\widetilde{\operatorname{mex}}}(r,n) q^n
	\end{align}
	Since $ \operatorname{mex}_{2}(\lambda)=2r-1,$ the compulsory odd parts are $ 1,3,5,\ldots,2r-3$ and their sum is $$1+3+\cdots+(2r-3)=(r-1)^2.$$
	Hence, these compulsory parts contribute
	$q^{(r-1)^2} \prod_{j=1}^{r-2}
	\frac{1}{1-q^{2j+1}}.$ The part
	$2r-1$
	is excluded, while every larger odd part may occur with arbitrary
	multiplicity. Thus,
	\begin{align}\label{eq100}
		\sum_{n=0}^{\infty}
		p_{2}^{\widetilde{\operatorname{mex}}}(r,n) q^n
		=
		q^{(r-1)^2} \prod_{j=1}^{r-2}
		\frac{1}{1-q^{2j+1}}
		\prod_{j=r}^{\infty}
		\frac{1}{1-q^{2j+1}} = \frac{q^{(r-1)^2}(1-q^{2r-1})} {(q;q^2)_\infty}.
	\end{align}
Combining \eqref{eq099} and \eqref{eq100} together, we get 
\begin{align}\label{eq100a}
	G_{2}(z,q)= \sum_{r=1}^{\infty} z^{2r-1} \sum_{n=0}^{\infty}
	p_{2}^{\widetilde{\operatorname{mex}}}(r,n) q^n = \sum_{r=1}^{\infty} z^{2r-1} \frac{q^{(r-1)^2}(1-q^{2r-1})} {(q;q^2)_\infty}. 
\end{align}
Differentiating $ 	G_{2}(z,q)$ with respect to $z$ and putting $z=1, $ we get
\begin{align} &\label{eq101}
		\sum_{n\geq0} \sigma_{2}^{\mex}(n)q^n = \left. \frac{\partial}{\partial z}G_{2}(z,q) \right|_{z=1} \\
		=& \sum_{r=1}^{\infty} (2r-1) \frac{ q^{(r-1)^2}(1-q^{2r-1}) } {(q;q^2)_\infty}
		= \frac{1}{{(q;q^2)_\infty}} \sum_{r=1}^{\infty} (2r-1) q^{(r-1)^2}(1-q^{2r-1})\notag\\
		=& \frac{1}{{(q;q^2)_\infty}} \Bigg[ \sum_{r=1}^{\infty} (2r-1) q^{(r-1)^2} - \sum_{r=1}^{\infty} (2r-1) q^{(r-1)^2}q^{2r-1}\Bigg]\notag  \\
		=& \frac{1}{{(q;q^2)_\infty}} \Bigg[ \sum_{t=0}^{\infty} (2t+1) q^{t^2} - \sum_{t=0}^{\infty} (2t+1) q^{(t+1)^2}\Bigg]\notag \\
		=& \frac{1}{{(q;q^2)_\infty}} \Bigg[ \sum_{t=0}^{\infty} (2t+1) q^{t^2} - \sum_{t=1}^{\infty} (2t-1) q^{t^2}\Bigg] \notag \\
		=& \frac{1}{{(q;q^2)_\infty}} \Bigg[1+ 2 \sum_{t=1}^{\infty}  q^{t^2}\Bigg] =  \frac{ \sum_{t\ \in \mathbb{Z}} q^{t^2}}{(q;q^2)_{\infty}}  = \frac{(q^2;q^2)^6_{\infty}}{(q;q)^3_{\infty} (q^4;q^4)^2_{\infty}}. \label{eq101a}
	\end{align}
\end{proof}
\subsection*{Proof of Theorem \ref{sigmamex3}}
In order to prove Theorem \ref{sigmamex3}, we need the following proposition.
	\begin{prop}\label{computationofA_3}
	For any positive integer $m,  \quad \hbox{with} \quad  3\nmid m,$ let 
	\begin{align}\label{equationA_m}
	 A_3(m)=\sum_{\substack{1\le j<m\\3\nmid j}}j, 
	\end{align}
	 Then
	\begin{align}
		A_{3r+1}=3r^2, \qquad  	A_{3r+2}=3r^2+3r+1 \quad \hbox{for some non-negative integer} \quad  r.
	\end{align}  
		\end{prop}
	\begin{proof}
		Now we compute the quantity  $ A_3(m).$  
	\paragraph{Case I}
	Suppose that $m=3r+1.$
	Then
	\begin{align*}
		A_3(3r+1)=\sum_{k=0}^{r-1}(3k+1)+\sum_{k=0}^{r-1}(3k+2).
	\end{align*}
	Since
	\begin{align*} \sum_{k=0}^{r-1}(3k+1) = \frac{r(3r-1)}{2} \qquad \hbox{and} \qquad 
		\sum_{k=0}^{r-1}(3k+2)
		=\frac{r(3r+1)}{2},
	\end{align*}
	we obtain
	$A_{3r+1}=3r^2.$
	\paragraph{Case II}
	Suppose that $m=3r+2.$
	Then
	\begin{align*}
		A_{3r+2}=A_{3r+1}+(3r+1)=3r^2+(3r+1).
	\end{align*}
\end{proof}
\begin{proof}[Proof of Theorem \ref{sigmamex3}]
From Proposition \ref{bivarategeneratingfunction}, we have
	\begin{align}\label{eq115}
		G_3(z,q)= & \sum_{n\ge0}
			\sum_{\substack{m\ge1\\3\nmid m}}
			p^{\mathrm{mex}}_{3}(m,n)
			z^m q^n =\sum_{\substack{m\ge1\\3\nmid m}} \frac{ z^{m} q^{A_3(m)}(1-q^m)}{(q;q^3)_{\infty} (q^2;q^3)_{\infty} }.
		\end{align}
Using Proposition \ref{computationofA_3} in \eqref{eq115}, we get
\begin{align}\label{eq115a}
 G_3(z,q)=  \frac{1}{(q;q^3)_\infty (q^2;q^3)_\infty }
 \Bigg[
 \sum_{r\ge0}
 z^{3r+1}q^{3r^2}\left(1-q^{3r+1}\right)
 +\sum_{r\ge0}z^{3r+2}q^{3r^2+3r+1}\left(1-q^{3r+2}\right)
 \Bigg].
\end{align}
	Differentiating \eqref{eq115a} with respect to $z$ and putting $z=1$, we obtain
	\begin{align}
		\begin{aligned}
			\frac{\partial }{\partial z} ( G_3(z,q)) \Bigg|_{z=1} &=
			\sum_{n\ge0}
			\sigma^{\mathrm{mex}}_{3}(n)q^n :=\sum_{\substack{m\ge1\\3\nmid m}}
			m\,p^{\mathrm{mex}}_{3}(m,n).\\
			&= \frac{1}{(q;q^3)_\infty (q^2;q^3)_\infty }
			\Bigg[
			\sum_{r\ge0}
			(3r+1)q^{3r^2}\left(1-q^{3r+1}\right)\\
			\qquad\qquad
			&+\sum_{r\ge0}(3r+2) q^{3r^2+3r+1}\left(1-q^{3r+2}\right)
			\Bigg]
			= \frac{N_3(q)}{(q;q^3)_\infty (q^2;q^3)_\infty },
		\end{aligned}
	\end{align}
where \begin{align}
		N_3(q) &= \sum_{r\ge0}
		(3r+1)q^{3r^2}\left(1-q^{3r+1}\right)
		+\sum_{r\ge0}(3r+2) q^{3r^2+3r+1}\left(1-q^{3r+2}\right) \\ \notag 
		&= \sum_{r\ge0}
		(3r+1)q^{3r^2} + \sum_{r\ge0}
		q^{3r^2 +3r+1}- \sum_{r\ge0}
		(3r+2)q^{3(r+1)^2} \\ \notag 
		&= \sum_{r\ge0}
		(3r+1)q^{3r^2} -\sum_{r\ge0}
		(3r+2)q^{3(r+1)^2} + \sum_{r\ge0}
		q^{3r^2 +3r+1}\\ \notag 
		&= \sum_{r\ge0}
		(3r+1)q^{3r^2} -\sum_{r\ge1}
		(3r-1)q^{3r^2} + \sum_{r\ge0}
		q^{3r^2 +3r+1} \\ \label{eq108} 
		&= 1+ 2 \sum_{r\ge1}
		q^{3r^2} +  \sum_{r\ge0}
		q^{3r^2 +3r+1}  = \sum_{r \in \mathbb{Z} } q^{3r^2} + \sum_{r\ge0}
		q^{3r^2 +3r+1}
	\end{align}
Note that, for any real number $r,$ \begin{align}\label{algebra_identity}
		(-r-1)^2+ (-r-1)
		= r^2+r.\end{align} 
	Using \eqref{algebra_identity}, we get
	\begin{align}\label{eq109}
		\sum_{r\ge0}
		q^{3r^2 +3r+1}= \frac{1}{2} \sum_{r \in \mathbb{Z}} q^{3r^2 +3r+1} = \frac{q}{2} \sum_{ r= - \infty}^{\infty}  \left( q^3\right)^{r^2} \left( q^3\right)^{r}.
	\end{align}
		Applying $Q=q^3$ and $z=q^3$  in Lemma \ref{lemmaJacobitripleproduct}, we have
	\begin{align}\label{eq110}
		\sum_{r=-\infty}^{\infty} q^{3r^2+3r}= \left(-q^6; q^6\right)_{\infty}\left(-1; q^6\right)_{\infty} \left(q^6; q^6\right)_{\infty} = 2 \left(-q^6; q^6\right)^2_{\infty}  \left(q^6; q^6\right)_{\infty},
	\end{align}
	since $\left(-1; q^6\right)_{\infty} = 2 \left(-q^6; q^6\right)_{\infty}.$
	Further observe that \begin{align}\label{eq111}
		\left(-q^6; q^6\right)_{\infty} = \frac{\left(q^{12}; q^{12}\right)_{\infty}}{\left(q^6; q^6\right)_{\infty}}.
	\end{align} Using \eqref{eq110} and \eqref{eq111} in \eqref{eq109}, we get
	\begin{align}
		\sum_{r\ge0}
		q^{3r^2 +3r+1}= q \frac{\left(q^{12}; q^{12}\right)^2_{\infty}}{\left(q^6; q^6\right)_{\infty}}.
	\end{align}
	Also, Using Lemma \ref{lem1}, we have
	\begin{align}
		\sum_{r \in \mathbb{Z} } q^{3r^2} = \frac{\left(q^6; q^6\right)^5_{\infty} }{\left(q^3; q^3\right)^2_{\infty} \left(q^{12}; q^{12}\right)^2_{\infty}}
	\end{align}
	Thus \begin{align}\label{expressionsigma3}
		\sum_{n\ge0}
		\sigma^{\mathrm{mex}}_{3}(n)q^n & =  \frac{1} {\left(q; q^3\right)_{\infty} \left(q^2; q^3\right)_{\infty}} \Bigg[ \frac{\left(q^6; q^6\right)^5_{\infty} }{\left(q^3; q^3\right)^2_{\infty} \left(q^{12}; q^{12}\right)^2_{\infty}} + q \frac{\left(q^{12}; q^{12}\right)^2_{\infty}}{\left(q^6; q^6\right)_{\infty}}.  \Bigg] \\ \notag 
		&=\frac{\left(q^3; q^3\right)_{\infty}} {\left(q; q\right)_{\infty} } \Bigg[ \frac{\left(q^6; q^6\right)^5_{\infty} }{\left(q^3; q^3\right)^2_{\infty} \left(q^{12}; q^{12}\right)^2_{\infty}} + q \frac{\left(q^{12}; q^{12}\right)^2_{\infty}}{\left(q^6; q^6\right)_{\infty}}.  \Bigg] \\ \notag
		&=  \frac{\left(q^6; q^6\right)^5_{\infty} }{ \left(q; q\right)_{\infty} \left(q^3; q^3\right)_{\infty} \left(q^{12}; q^{12}\right)^2_{\infty}} + q \frac{ \left(q^{12}; q^{12}\right)^2_{\infty} \left(q^{3}; q^{3}\right)_{\infty} }{\left(q;q\right)_{\infty} \left(q^6; q^6\right)_{\infty}}.  
	\end{align}
\end{proof}	
\begin{proof}[Proof of Theorem \ref{monotonicitytheorem}]
		Consider the map $ \phi:\mathcal{P}_{\ell}(n) \longrightarrow \mathcal{P}_{\ell}(n+1) $ by $ \phi(\lambda)=\lambda\cup\{1\}. $
	Since $1$ is never divisible by $\ell$, $\phi(\lambda)$ is still an $\ell$-regular partition. Moreover, $ |\phi(\lambda)| = |\lambda|+1 = n+1. $ The map is injective, because if $ \phi(\lambda_1)=\phi(\lambda_2), $ then removing one copy of $1$ from both partitions gives $ \lambda_1=\lambda_2. $
	Now, we consider two cases. Firstly, suppose that $1\in\lambda,$ then $ \mex_{\ell}(\lambda)>1. $ Adding another $1$ does not change whether any other part occurs. Hence $ \mex_{\ell}(\phi(\lambda)) = \mex_{\ell}(\lambda). $ Secondly, suppose that $1\notin\lambda,$ then $\mex_{\ell}(\lambda)=1.$ After adding $1$, $ 1\in\phi(\lambda).$ 
	Therefore, the mex of $\phi(\lambda)$ is the smallest positive integer other than $1$ which is not divisible by $\ell$ and does not occur in $\lambda$. Thus, $ \mex_{\ell}(\phi(\lambda)) > 1 = \mex_{\ell}(\lambda). $ 
	Therefore, in both cases, $\mex_{\ell}(\phi(\lambda)) \geq \mex_{\ell}(\lambda).$ Since $\phi$ is injective, we have
	\begin{align*} \sigma_{\ell}^{\mex}(n+1) &= \sum_{\mu\in\mathcal{P}_{\ell}(n+1)} \mex_{\ell}(\mu) \geq \sum_{\lambda\in\mathcal{P}_{\ell}(n)} \mex_{\ell}(\phi(\lambda)) \geq \sum_{\lambda\in\mathcal{P}_{\ell}(n)} \mex_{\ell}(\lambda) = \sigma_{\ell}^{\mex}(n). \end{align*}
	Hence $  \sigma_{\ell}^{\mex}(n+1) \geq \sigma_{\ell}^{\mex}(n) \quad (n\geq 0). $ Therefore, the sequence $ \left\{\sigma_{\ell}^{\mex}(n)\right\}_{n\geq 0} $ is weakly increasing for every $\ell\geq 2$.
	
As defined earlier, we have	 $ p_{\ell}(n) = \left|\mathcal{P}_{\ell}(n)\right|. $ 
	 Every partition $\lambda\in\mathcal{P}_{\ell}(n)$ with $1\notin\lambda$ satisfies $ \mex_{\ell}(\lambda)=1. $
	After adding a part $1$, its mex becomes strictly larger, i.e.  $ \mex_{\ell}(\phi(\lambda)) > 1 = \mex_{\ell}(\lambda). $ 
	If $ 1 \in \lambda,$ then $ \mex_{\ell}(\phi(\lambda))  = \mex_{\ell}(\lambda). $
	Hence
	\begin{align}\label{eqquantitative1}
		\sigma_{\ell}^{\mex}(n+1) - \sigma_{\ell}^{\mex}(n) \geq \#\left\{ \lambda\in\mathcal{P}_{\ell}(n):1\notin\lambda \right\}. 
	\end{align}
	Now consider the partitions in $\mathcal{P}_{\ell}(n)$ which contain at least one part $1$. 
	Removing one occurrence of $1$ gives a bijection with $\mathcal{P}_{\ell}(n-1)$. 
	Therefore,
	\begin{align}\label{eqquantitative2}
		\#\left\{ \lambda\in\mathcal{P}_{\ell}(n):1\notin\lambda \right\} = p_{\ell}(n)-p_{\ell}(n-1).
	\end{align}
	Comparing \eqref{eqquantitative1} and \eqref{eqquantitative2}, we get 
	\begin{align*}
		\sigma_{\ell}^{\mex}(n+1) - \sigma_{\ell}^{\mex}(n) \geq p_{\ell}(n)-p_{\ell}(n-1)
	\end{align*}
\end{proof}

%

\begin{proof}[Proof of Theorem \ref{asymptotictheoremofsigma2}]
Let \[C(q)= \sum_{n\geq 0} \sigma_{2}^{\mex}(n) q^n 	=
	\frac{(q^2;q^2)_\infty^6}
	{(q;q)_\infty^3(q^4;q^4)_\infty^2}.
	\]
Applying the transformation formula for the Dedekind's eta function \cite{Koblitz} , we have
\begin{align}\label{Dedekindtransformation}
	(e^{-y};e^{-y})_\infty
	\sim
	\sqrt{\frac{2\pi}{y}}
	\exp\left(-\frac{\pi^2}{6y}\right),
	\qquad y\to0^+.
\end{align}
	Thus
	\[
	(e^{-2y};e^{-2y})_\infty^6
	\sim
	\left(\sqrt{\frac{\pi}{y}}\right)^6
	\exp\left(-\frac{\pi^2}{2y}\right)
	=
	\frac{\pi^3}{y^3}
	e^{-\frac{\pi^2}{2y}},
	\]
	while
	\[
\left(e^{-y};e^{-y}\right)_\infty^3
	\sim
	\left(\sqrt{\frac{2\pi}{y}}\right)^3
	e^{-\frac{\pi^2}{2y}},
	\]
	and
	\[
	(e^{-4y};e^{-4y})_\infty^2
	\sim
	\left(\sqrt{\frac{\pi}{2y}}\right)^2
	e^{-\frac{\pi^2}{12y}}
	=
	\frac{\pi}{2y}
	e^{-\frac{\pi^2}{12y}}.
	\]	Therefore,
	\begin{align*}
		C(e^{-y}) \sim  \sqrt{\frac{\pi}{2}}\,
		y^{-1/2}
		e^{\frac{\pi^2}{12y}}.
	\end{align*}
From Theorem \ref{monotonicitytheorem}, we know that the sequence $\{ \sigma_{2}^{\mex}(n) \}$ is weakly increasing. Thus, by applying Lemma \ref{LemmaIngham} with parameter 

	$\beta=\sqrt{\frac{\pi}{2}},\qquad \gamma=-\frac12, \qquad \mu=\frac{\pi^2}{12},$ we obtain 
		\begin{align*}
	\sigma_{2}^{\mex}(n)
	\sim
	\frac{1}{2\sqrt{2n}}\,
	e^{\pi\sqrt{n/3}},
	\qquad \hbox{as} \quad  n\to\infty.
	\end{align*}
\end{proof}
\begin{proof}[Proof of Theorem \ref{asymptotictheoremofsigma3}]
 From \eqref{expressionsigma3},	we have
	\[
	\sum_{n\geq0}
	\sigma^{\mathrm{mex}}_{3}(n)q^n
	=\frac{\left(q^6;q^6\right)^5_{\infty}}
	{\left(q;q\right)_{\infty}
		\left(q^3;q^3\right)_{\infty}
		\left(q^{12};q^{12}\right)^2_{\infty}}
	+
	q\frac{
		\left(q^{12};q^{12}\right)^2_{\infty}
		\left(q^3;q^3\right)_{\infty}}
	{\left(q;q\right)_{\infty}
		\left(q^6;q^6\right)_{\infty}}.
\]
	
%
	
Let
	\begin{align*}
		T_1(q):= \frac{(q^6;q^6)_\infty^5} {(q;q)_\infty(q^3;q^3)_\infty(q^{12};q^{12})_\infty^2}.
	\end{align*}
	Using \eqref{Dedekindtransformation}, we get
	\[
	(e^{-6y};e^{-6y})_\infty^5
	\sim
	\left(\sqrt{\frac{\pi}{3y}}\right)^5
	e^{-\frac{5\pi^2}{36y}},
	\]
	and
	\[
	(e^{-y};e^{-y})_\infty(e^{-3y};e^{-3y})_\infty(e^{-12y};e^{-12y})_\infty^2
	\sim
	\sqrt{\frac{2\pi}{y}}
	\sqrt{\frac{2\pi}{3y}}
	\left(\sqrt{\frac{\pi}{6y}}\right)^2
	e^{-\frac{\pi^2}{6y}
		-\frac{\pi^2}{18y}
		-\frac{\pi^2}{36y}}.
	\]	Therefore,
	\begin{align*}
T_1(e^{-y})
	\sim
	\frac{\sqrt{\pi}}{3}\,
	y^{-1/2}
	e^{\frac{\pi^2}{9y}}.
\end{align*}
	Similarly, let
	\begin{align*}
	T_2(q):=
	q\frac{
		(q^{12};q^{12})_\infty^2(q^3;q^3)_\infty}
	{(q;q)_\infty(q^6;q^6)_\infty}.
\end{align*}
	
As $y\to0^+,$ we get $q=e^{-y}\to1$ and hence $q\sim1$. Thus
	\[
	T_2(e^{-y})
	\sim
	\frac{\sqrt{\pi}}{6}\,
	y^{-1/2}
	e^{\frac{\pi^2}{9y}}.
	\] Note that both terms have the same exponential growth. Consequently,
	\begin{align*}
T_1(e^{-y}) + T_2(e^{-y}) =\sum_{n\geq0}
		\sigma^{\mathrm{mex}}_3(n)e^{-ny}
		\sim
		\left(
		\frac{\sqrt{\pi}}{3}
		+
		\frac{\sqrt{\pi}}{6}
		\right)
		y^{-1/2}
		e^{\frac{\pi^2}{9y}}=
		\frac{\sqrt{\pi}}{2}
		y^{-1/2}
		e^{\frac{\pi^2}{9y}}.
	\end{align*}
	Applying the Lemma \ref{LemmaIngham} with parameter
	$ \beta=\frac{\sqrt{\pi}}{2},
	\quad
	\gamma=-\frac12,
	\quad
	\mu=\frac{\pi^2}{9},$ we get
	
%
%
\begin{align*}
\sigma^{\mathrm{mex}}_3(n)
	\sim
	\frac{1}{4\sqrt{n}}\,
	e^{\frac{2\pi}{3}\sqrt{n}},
	\qquad \hbox{as} \qquad n\to\infty.
\end{align*}
	\end{proof} 
\subsection{Proof of Theorem \ref{distributionofsigma2} }
In oredr to prove Theorem \ref{distributionofsigma2}, we need the following parity result for $\sigma_2^{\mathrm{mex}}(n).$
\begin{prop}\label{pairtyresultsigma2} We have
			\begin{align*}
		\sigma_2^{\mathrm{mex}}(n)\equiv 1\pmod{2}
		\quad \hbox{ if and only if } \quad n=0 \quad\text{or}\quad
		n=\frac{k(3k-1)}{2}
		\quad\text{or}\quad
		n=\frac{k(3k+1)}{2}
			\end{align*} for some positive integer $k.$
	\end{prop}
\begin{proof}
	From equation \eqref{eq101a}, we have 
	\begin{align}\label{eq101b}
\sum_{n=0}^{\infty} \sigma_2^{\mathrm{mex}}(n) q^n =	\frac{(q^2;q^2)^6_{\infty}}{(q;q)^3_{\infty} (q^4;q^4)^2_{\infty}}.
	\end{align}
 By using Lemma \ref{lemmabinomial}, we get
\begin{align}\label{eq101c}
	\sum_{n=0}^{\infty} \sigma_2^{\mathrm{mex}}(n) q^n \equiv 	(q;q)_{\infty} \pmod 2.
\end{align}
From Lemma \ref{pentagonaltheorems}, taking modulo $2,$ we obtain
\begin{align}\label{eq101v}
(q;q)_{\infty} \equiv
1+\sum_{k\geq 1}
\left(
q^{\frac{k(3k-1)}{2}}
+
q^{\frac{k(3k+1)}{2}}
\right)
\pmod{2}.
\end{align}
From \eqref{eq101v} and \eqref{eq101c}, we obtain
\begin{align*}
\sigma_2^{\mathrm{mex}}(n)\equiv 1\pmod{2}
\quad\Longleftrightarrow\quad
n=\frac{k(3k-1)}{2}
\quad\text{or}\quad
n=\frac{k(3k+1)}{2}
\end{align*}
for some integer \(k\geq 1\).
	\end{proof}

\begin{proof}[Proof of Theorem \ref{distributionofsigma2}]
 For a real number $X,$	we define
	\[
	A(X):=\#\{n\leq X: \sigma_2^{\mathrm{mex}}(n) \text{ is odd}\},
	\]
	and
	\[
	B(X):=\#\{n\leq X:\sigma_2^{\mathrm{mex}}(n)\text{ is even}\}.
	\]	
Note that for any real number $k,$
	\[
	\frac{k(3k-1)}{2}\leq X \iff k\leq \frac{1+\sqrt{1+24X}}{6}. 
	\] Similarly, 
	\[
	\frac{k(3k+1)}{2}\leq X \iff k\leq \frac{-1+\sqrt{1+24X}}{6}.
	\]
	Using Proposition \ref{pairtyresultsigma2}, we get
	\begin{align}
	A(X) &= \left\lfloor \frac{1+\sqrt{1+24X}}{6}\right \rfloor
	+
	\left\lfloor\frac{-1+\sqrt{1+24X}}{6}\right\rfloor\\ \notag 
	& = \frac{1+\sqrt{1+24X}}{6} + \frac{-1+\sqrt{1+24X}}{6} + \mathcal{O}(1)= \frac{\sqrt{1+24X}}{3} + \mathcal{O}(1).
	\end{align} and
	\[
	  	B(X)=\lfloor X\rfloor-A(X) = X- A(x)+ \mathcal{O}(1).
	\]
	
For any real number $u$, we have
\begin{align}\label{binomialexpansion}
(1+u)^{1/2}
=
1+\frac{1}{2}u-\frac{1}{8}u^2
+\frac{1}{16}u^3
-\frac{5}{128}u^4+\cdots,
\end{align}
	For large \(x\), using \eqref{binomialexpansion}, we get
	\begin{align}
\sqrt{x+a} =
	\sqrt{x}\left(1+\frac{a}{x}\right)^{1/2} = \sqrt{x} +\frac{a}{2\sqrt{x}}
	-\frac{a^2}{8x^{3/2}}
	+\frac{a^3}{16x^{5/2}}
	-\frac{5a^4}{128x^{7/2}}
	+\cdots.
\end{align}
For large value of $x,$
	\[ \sqrt{x+a} =
	\sqrt{x}
	+\frac{a}{2\sqrt{x}}
	+ \mathcal{O}\left(x^{-3/2}\right).
	\]
In particular	
	\[
	\frac{\sqrt{24X+1}}{3}
	=
	\frac{\sqrt{24X}}{3}
	+\frac{1}{6\sqrt{24X}}
	+\mathcal{O} \left(X^{-3/2}\right).
	\]Thus
	\[
	A(X)=\frac{2\sqrt{6}}{3}\sqrt{X}+\mathcal{O}(1). \quad 
	\]
	Therefore,
	\[
	B(X)
	=
	X-\frac{2\sqrt{6}}{3}\sqrt{X}+\mathcal{O}(1).
	\]
	\end{proof}

\subsection{Proof of Theorem \ref{distributionofsigma3}}
In order to prove Theorem \ref{distributionofsigma3}, we need the following lemma.

	\begin{lemma}\label{disributionfirstlemmasigma3}
	Let $n$ be a positive integer such that  $n\equiv2\pmod{3}.$ Then we have
	\begin{align*}
		\sum_{s=0}^{\infty}
		\sigma_3^{\mathrm{mex}}
		\left(n-\frac{s(3s-1)}{2}\right)
		+
		\sum_{s=1}^{\infty}
		\sigma_3^{\mathrm{mex}}
		\left(n-\frac{s(3s+1)}{2}\right)
		\equiv0\pmod{2}.
	\end{align*}
\end{lemma}
\begin{proof}
 Let \begin{align}\label{eq300}
  F(q):=\sum_{n\geq 0}\sigma_3^{\mathrm{mex}}(n)q^n =  \frac{\left(q^6; q^6\right)^5_{\infty} }{ \left(q; q\right)_{\infty} \left(q^3; q^3\right)_{\infty} \left(q^{12}; q^{12}\right)^2_{\infty}} + q \frac{ \left(q^{12}; q^{12}\right)^2_{\infty} \left(q^{3}; q^{3}\right)_{\infty} }{\left(q;q\right)_{\infty} \left(q^6; q^6\right)_{\infty}}.  
\end{align}
	Taking modulo \(2\) of \eqref{eq300}, we get
	\begin{align*}
	F(q)\equiv \frac{f_3}{f_1}\left(1+qf_3^6\right)\pmod{2} \iff f_1F(q)\equiv f_3+qf_3^7\pmod{2}.
	\end{align*}
Thus,
	\begin{align}\label{eq302}
	(q;q)_\infty
	\sum_{n\geq 0}\sigma_3^{\mathrm{mex}}(n)q^n
	\equiv f_3+qf_3^7\pmod{2}.
\end{align}


From Lemma \ref{pentagonaltheorems} and taking	modulo \(2\) on it, we get
	\begin{align}\label{eq304}
	(q;q)_\infty \equiv
	1+\sum_{s=1}^{\infty}q^{s(3s-1)/2}
	+\sum_{s=1}^{\infty}q^{s(3s+1)/2}
	\pmod{2}.
\end{align}
	
	Using \eqref{eq304} in \eqref{eq302}, we obtain
	\begin{align}\label{eq305}
	\left( \sum_{s=0}^{\infty}q^{s(3s-1)/2}
	+\sum_{s=1}^{\infty}q^{s(3s+1)/2}
	\right)
	\left(
	\sum_{n\geq0}\sigma_3^{\mathrm{mex}}(n)q^n
	\right)
	\equiv
	f_3+qf_3^7
	\pmod{2}.
\end{align}
Since $qf_3^7$ contains only powers $q^{3j+1}$ and comparing the coefficients of \(q^n\) of \eqref{eq305}, we obtain
	\begin{align}\label{eq307}
	\sum_{s=0}^{\infty}
	\sigma_3^{\mathrm{mex}}
	\left(n-\frac{s(3s-1)}{2}\right)
	+
	\sum_{s=1}^{\infty}
	\sigma_3^{\mathrm{mex}}
	\left(n-\frac{s(3s+1)}{2}\right)
	\equiv0\pmod{2},
	\qquad n\equiv2\pmod{3}.
	\end{align}
		\end{proof} 
	
\begin{lemma}\label{existancelemma}
	For every integer \(r\geq2\) satisfying $ r\equiv2\pmod{3}, $
	there exists an integer $ n\in \left[ 2r-1,\frac{r(3r-1)}{2} \right] $
	such that $ \sigma_3^{\mathrm{mex}}(n)\equiv1\pmod{2}.$ In particular, there are infinitely many integers $n$ for which $ 	\sigma_3^{\mathrm{mex}}(n) $ is an odd integer.
\end{lemma} 
\begin{proof} 
	 On contrary suppose that $ \sigma_3^{\mathrm{mex}}(n) $
	is even for every $ 2r-1\leq n\leq\frac{r(3r-1)}{2}. $
	We define
	\begin{align}
S(s):=\frac{r(3r-1)}{2} -\frac{s(3s-1)}{2} \quad \hbox{and} \quad  R(s):=\frac{r(3r-1)}{2}
-\frac{s(3s+1)}{2}.
\end{align}
Observe that if $ r\equiv2\pmod{3}, $ then 	$\frac{r(3r-1)}{2} \equiv r\pmod{3}.$
Therefore, \eqref{eq307} can be rewritten as 
	\begin{align}\label{eq308}
\sum_{s=0}^{\infty}
	\sigma_3^{\mathrm{mex}}(S(s))
	+
	\sum_{s=1}^{\infty}
	\sigma_3^{\mathrm{mex}}(R(s))
	\equiv0\pmod{2}, \quad r\equiv2\pmod{3}.
\end{align}
Notice that for \(s\geq r+1\), $ S(s)<0, $
	and for \(s\geq r\), $ R(s)<0.$ We use the standard convention
$ \sigma_3^{\mathrm{mex}}(n)=0 \qquad (n<0).$
	Thus \eqref{eq308} reduces to
		\begin{align}\label{eq310}
		1+	\sum_{s=0}^{r-1} \sigma_3^{\mathrm{mex}}(S(s)) + \sum_{s=1}^{r-1}
		\sigma_3^{\mathrm{mex}}(R(s)) \equiv0\pmod{2},
	\end{align} since
	$ \sigma_3^{\mathrm{mex}}(S(r)) = \sigma_3^{\mathrm{mex}}(0) =1.$
Now we show that $S(s),R(s)$ lie in the interval $\left[2r-1,\frac{r(3r-1)}{2}\right]$  for \(1\leq s\leq r-1.\) Notice that $ S(0)=\frac{r(3r-1)}{2}, \quad S(r-1) = 3r-2.$ Note that \(S(s)\) decreases in \(s\) and the smallest value for
	\(0\leq s\leq r-1\) is \(S(r-1)\). Since $ 3r-2\geq2r-1 \quad (r\geq1),$
	we have
$ 2r-1\leq S(s)\leq\frac{r(3r-1)}{2}, \quad \hbox{ for} \quad 0\leq s\leq r-1. $
Similarly, for \(1\leq s\leq r-1\), again \(R(s)\) decreases with \(s\), and $ R(r-1) = 2r-1. $ Thus $ 2r-1\leq R(s)\leq\frac{r(3r-1)}{2}.$ Therefore, $S(s),R(s)$ lie in the interval $\left[2r-1,\frac{r(3r-1)}{2}\right].$ By our assumption, the left-hand side of \eqref{eq310} is
		$1+0+0\equiv1\pmod{2}, $  contradicting \eqref{eq310}, which says that it is \(0\pmod{2}\). Therefore, at least one integer \(n\) from $ n\in \left[ 2r-1,\frac{r(3r-1)}{2} \right] $
		such that $ \sigma_3^{\mathrm{mex}}(n)\equiv1\pmod{2}.$
	\end{proof}
\subsection*{Proof of Theorem \ref{distributionofsigma3}}	
	Suppose
	$ a_k\equiv2\pmod{3},$ for some  positive integer $k.$ Then $ a_{k+1} = \frac{a_k(3a_k-1)}{2}  \equiv2\pmod{3}.$
	Thus, $ a_k\equiv2\pmod{3} \quad\text{for every }k\geq0.$ In particular, we choose $a_0=2, \quad a_1= 5, \quad  a_2= 35, $ and so on. By applying the Lemma \ref{existancelemma} with \(r=a_k\), there exists some integer \(n_k\) such that
	\begin{align}\label{eq311}
 a_k \leq  2a_k-1 \leq n_k \leq \frac{a_k(3a_k-1)}{2} = a_{k+1} \quad \hbox{and} \quad \sigma_3^{\mathrm{mex}}(n_k) \equiv1\pmod{2}.
\end{align}
Let \(\nu\) be the largest integer such that $ a_\nu\leq X.$ Then
$a_\nu\leq X<a_{\nu+1}.$ From \eqref{eq311}, we have $ n_k\geq2a_k-1.$
	Since $a_{k+1} = \frac{a_k(3a_k-1)}{2},$ we have, for \(a_k\geq2\), $a_{k+1}>2a_k.$ Thus the interval $[2a_k-1,a_{k+1}]$ is disjoint with the analogous interval corresponding to \(k+2.\)  Now the number of $n_k$ for which $\sigma_3^{\mathrm{mex}}(n_k)$ is odd is at least ${\nu}.$ Now for all $k \geq 0,$ we obtain
	$$ a_{k+1}=\frac{a_k(3a_k-1)}{2} < \frac{3}{2}a_k^2 < 2 a_k^2. $$
On iteration, we get $ a_k<(2a_0)^{2^k} < 5^{2^k} \implies \log\log a_k < k\log2+\log\log 5.$ Since $ a_\nu\leq X<a_{\nu+1},$ we have
$ \log\log X<\log\log a_{\nu+1} <
\nu\log2+C$ where $ C=\log\log 5 + \log 2 $ is a fixed constant.
Therefore, for sufficiently large \(X\),$ \nu\gg\log\log X. $ \hfill $\square$

\section{Conclusion}

In this paper, we have introduced and studied the minimal excludant statistic for $\ell$-regular partitions, extending the classical work of Andrews and Newman to this more general setting. Further, we obtained explicit eta-quotient representations for its generating function when $\ell = 2$ and $\ell = 3$, revealing an elegant connection between this combinatorial statistic and the theory of modular forms. We also obtained asymptotic formula for $\sigma_2^{\mathrm{mex}}(n)$ and $\sigma_3^{\mathrm{mex}}(n).$ Further we obtained that the parity result for $\sigma_2^{\mathrm{mex}}(n)$ is governed entirely by the generalized pentagonal numbers which helps us to find a density result for  $\sigma_2^{\mathrm{mex}}(n)$ while $\sigma_3^{\mathrm{mex}}(n)$ is shown to be odd infinitely often, with a quantitative $\log\log X$ lower bound obtained through a delicate dissection argument involving Euler's pentagonal number theorem modulo $2$. These results suggest that the minimal excludant statistic for $\ell$-regular partitions carries rich arithmetic structure worthy of further exploration. Natural directions for future work include extending the parity and density results to $\sigma_\ell^{\mathrm{mex}}(n)$ for $\ell \geq 4$, investigating congruences modulo primes other than $2$, and studying the mex statistic for other classes of restricted partitions.

\noindent{\bf Data availability statement:} There is no data associated to our manuscript.

\end{document}